\documentclass[11pt]{article}

\usepackage{amsfonts,amssymb,amsmath,amsthm,epsfig,euscript}
\usepackage{graphicx}
\usepackage{epstopdf}
\newtheorem{theorem}{Theorem}
\newtheorem{lemma}[theorem]{Lemma}

\long\def\symbolfootnote[#1]#2{\begingroup
\def\thefootnote{\fnsymbol{footnote}}\footnote[#1]{#2}\endgroup}

\newcommand{\fig}[2]{\begin{figure}[ht]
\centering\includegraphics[scale=.66]{#1}
\caption{#2}
\label{figure:#1}
\end{figure}}

\title{Ranking and unranking trees with given degree 
sequences}

\author{Jeffery B.  Remmel\\
\small Department of Mathematics\\[-0.8ex]
\small U.C.S.D., La Jolla, CA, 92093-0112\\[-0.8ex]
\small \texttt{jremmel@ucsd.edu}\\
\and
S. Gill Williamson\\
\small Department of Computer Science and Engineering\\[-0.8ex]
\small U.C.S.D., La Jolla, CA, 92093-0114\\[-0.8ex]
\small \texttt{gwilliamson@ucsd.edu}
}

\date{ %\small Submitted:   ; Accepted: \\
\small MR Subject Classifications: 05A15, 05C05, 05C20, 05C30}

\begin{document}
\maketitle

\begin{abstract}In this paper, we provide algorithms to rank and unrank
certain degree-restricted classes of Cayley trees.  Specifically, we
consider classes of
trees that have a given degree sequence or a given multiset of degrees. Using special properties of a bijection due to E\u gecio\u glu and
Remmel \cite{ER1},
we show that one can reduce the problem of ranking and unranking these
classes of
degree-restricted trees to corresponding problems of ranking and
unranking certain classes of set partitions.  If the underlying set of 
trees have $n$ vertices, then the largest ranks involved in each case are of 
order $n!$ so that it takes $O(nlog(n)$ bits just to write down the ranks. 
Our ranking and unranking algorithms  
for these degree-restricted classes are as efficient as can be expected 
since we show that 
they require $O(n^2log(n)$ bit operations if the underlying trees 
have $n$ vertices. 
%\footnotemark
%\keywords{Trees, Forests, Ranking, Unranking, Algorithms.}
\end{abstract}

\section{Introduction}
In computational combinatorics, it is important to be able to {\it
efficiently\/} rank,
unrank, and randomly generate (uniformly) basic classes of combinatorial
objects. A ranking
algorithm for a finite set $S$ is a bijection from $S$ to the set 
$\{0, \cdots,|S|-1\}$.  An unranking algorithm is the inverse of a 
ranking algorithm.
Ranking and unranking techniques are useful for storage and 
retrieval of elements of
$S$.  Uniform
random generation plays a role in Monte Carlo methods and in search
algorithms such as hill
climbing or genetic algorithms over classes of combinatorial objects.
Uniform random generation of objects is always possible if 
one has an unranking algorithm since one can generate, 
uniformly, an integer in $\{0, \cdots, |S|-1\}$ and unrank.

Given the set $V=\{1, \cdots, n\}$, we consider the set
$C_n$ of trees with vertex set $V$.  These trees are sometimes called
Cayley trees and can
be viewed as the set of spanning trees of the complete graph $K_n$.
Ranking and unranking
algorithms for the set $C_n$ have been described by many authors. Indeed,
efficient ranking
and unranking algorithms have been given for classes of trees and forests
that considerably
generalize the Cayley trees (e.g., [3], [4], [5], [6], [7]).

In this paper we consider more refined problem,namely, ranking and unranking subsets of $C_n$ with specified degree sequences or a specified multisets of degrees.  The resolution of this refined 
problem hinges on having a bijection with certain very special properties. Let
$\vec{C}_{n,1}$ be the set of directed trees on $V$ that are rooted at 1. That 
is, a directed tree $T \in \vec{C}_{n,1}$ has all its edges 
directed towards its root 1. We replace $C_n$ with the equivalent set
$\vec{C}_{n,1}$.  For any tree $T \in C_n$,
$\sum_{i=1}^n deg_T(i) = 2n-2$.
Let $\vec{s} = \langle s_1, \ldots, s_n \rangle$ be any sequence of
positive integers such
that $\sum_{i=1}^n s_i =2n-2$. Then we  let 
$\vec{C}_{n,\vec{s}} = \{T \in 
\vec{C}_{n,1}:
\langle deg_T(1),
\ldots, deg_T(n)\rangle = \vec{s}\}$. It will 
easily follow from our results in section 2 that 
\begin{equation}\label{count1}
|\vec{C}_{n,\vec{s}}| = \binom{n-2}{{s_1-1, \ldots,s_n-1}}.
\end{equation}
Similarly if $S =\{1^{\alpha_1}, \ldots,
(n-1)^{\alpha_{n-1}}\}$ is a multiset such that $\sum_{i=1}^{n-1} \alpha_i
\cdot i = 2n-2$ and $\sum_{i=1}^n \alpha_i = n$,
then we let
$\vec{C}_{n,S} = \{T \in \vec{C}_{n,1}: \{deg_T(1), \ldots, deg_T(n)\} = S \}$.  It is easy to see from (\ref{count1}) that 
\begin{equation}\label{count2}
|\vec{C}_{n,S}| = \binom{n}{{\alpha_1, \ldots, \alpha_n}}\binom{n-2}{{s_1-1, \ldots,s_n-1}}.
\end{equation}

In \cite{ER1}, E\u gecio\u glu and Remmel 
constructed a fundamental bijection $\Theta$
between the class
of functions ${\cal F}_n = \{f:\{2, \ldots , n-1\}\rightarrow \{1, \ldots,
n\}\}$ and the set $\vec{C}_{n,1}$. We shall exploit a key property of 
this bijection which is that for any vertex $i$, $1 +|f^{-1}(i)|$ equals 
that degree of $i$ in the tree $T= \Theta(f)$ when $\Theta(f)$ is considered 
as an undirected graph. This property allows us to reduce the problem 
or ranking and unranking trees in $\vec{C}_{n,\vec{s}}$ or $\vec{C}_{n,S}$ to 
the problem of ranking and unranking certain set partitions.  
We then use known techniques for ranking and unranking set partitions 
\cite{W}.

Note that both the $|\vec{C}_{n,\vec{s}}|$ and $|\vec{C}_{n,S}|$ can be as large as 
as the order of $n!$ so that the numbers involved in ranking and unranking will require on the order of $nlog(n)$ bits just to write down. 
We shall show that our algorithms for the ranking and unranking algorithms for
$\vec{C}_{n,\vec{s}}$ and $\vec{C}_{n,S}$ are as efficient as can be expected in that given a tree $T$, it will require at most $O(n^2)$ comparisons of 
numbers $y \leq n$ plus $O(n)$ operations of multiplication, division, addition, substraction and comparision on numbers 
$x < |\vec{C}_{n,\vec{s}}|$ ($x < |\vec{C}_{n,S}|$)
to find the rank of $T$ in $\vec{C}_{n,\vec{s}}$ ($\vec{C}_{n,S}$) so 
that it will take $O(n^2log(n))$ bit operations for ranking in 
either $\vec{C}_{n,\vec{s}}$ or $\vec{C}_{n,S}$.  Similarly, we will show that that the unranking algorithms $\vec{C}_{n,\vec{s}}$ or $\vec{C}_{n,S}$ will 
reqire  at most $O(n^2log(n))$ bit operations.

The outline of this paper is as follows. In Section~2, we describe
the bijection $\Theta:{\cal F}_n \rightarrow \vec{C}_{n,1}$ of \cite{ER1}
and discuss some
of its key properties.  In Section~3, we show that both $\Theta$ and
$\Theta^{-1}$
can be computed in linear time. This result allows us to reduce the problem of
efficiently ranking and unranking trees in $\vec{C}_{n,\vec{s}}$ or 
$\vec{C}_{n,S}$ 
to the
problem of efficiently ranking and unranking certain classes of set
partitions.  In
Section~4,  we  shall give ranking and unranking algorithms for the
classes of set
partitions corresponding the sets $\vec{C}_{n,\vec{s}}$ and $\vec{C}_{n,S}$.

\section{The $\Theta$ Bijection and its Properties}

In this section, we shall review the bijection 
$\Theta: {\cal F}_n \rightarrow \vec{C}_{n,1}$ due to 
E\u gecio\u glu and Remmel \cite{ER1} and give some of its properties. 

Let $[n] = \{1, 2, \ldots, n\}$. For each function 
$f:\{2, \ldots, n-1\} \rightarrow [n]$, we associate a directed graph $f$, 
$graph(f) = ([n],E)$ by setting $E = \{\langle i,f(i) \rangle :i = 2, \ldots, n-1\}$.   Following \cite{RW}, given  any  directed edge $(i,j)$
where $1 \leq i,j \leq n$, we define the weight of $(i,j)$, $W((i,j))$, by
\begin{equation}\label{wedge}
W((i,j)) =
\left\{ \begin{array}{ll}
p_is_j \  \mbox{if $i < j$}, \\
q_it_j \  \mbox{if $i \geq j$}
\end{array} \right.
\end{equation}
where $p_i, q_i, s_i, t_i$ are variables for $i = 1, \ldots, n$.
We shall call a directed edge $(i,j)$ a {\em descent edge} if $i \geq j$ and an
{\em ascent edge} if $i < j$.  We then define the weight of any digraph $G =
([n], E)$ by
\begin{equation}\label{wgraph}
W(G) = \prod_{(i,j) \in E} W((i,j)).
\end{equation}
A moment's thought will convince one that, in
general, the digraph corresponding to a function $f \in {\cal F}_n$
will consists of $2$ root-directed trees
rooted at vertices $1$ and $n$ respectively, with all
edges directed toward their roots, plus a number of
directed cycles of length $\geq 1$. For each vertex $v$ on a given cycle,
there is possibly a root-directed tree
attached to $v$  with $v$ as the root and all edges directed toward $v$.
Note the fact that there are trees rooted at vertices
$1$ and $n$ is due to the fact that these elements are not in the domain
of $f$.  Thus 
there can be no directed edges out of any of these vertices. We let the
weight of
$f$, $W(f)$, be the weight of the digraph $graph(f)$ associated with $f$.

To define the bijection $\Theta$, we first imagine
that the directed graph corresponding
to $f \in {\cal F}$ is drawn  so that
\begin{description}
\item[(a)] the trees rooted at $n$ and $1$ are drawn on the extreme left
and the extreme
right respectively with their edges directed upwards,

\item[(b)] the cycles are drawn so that their vertices form a directed path
on the line
between $n$ and $j$, with one back edge above the line, and  the
root-directed tree attached
to any vertex on a cycle is
drawn below the line between $n$ and $1$  with its edges directed upwards,

\item[(c)] each cycle $c_i$ is arranged so that its maximum
element $m_i$ is on the right, and 

\item[(d)] the cycles are arranged from left to right by decreasing maximal elements.
\end{description}
Figure~1 pictures a function $f$ drawn according to the rules  
(a)-(d) where $n =23$.

\fig{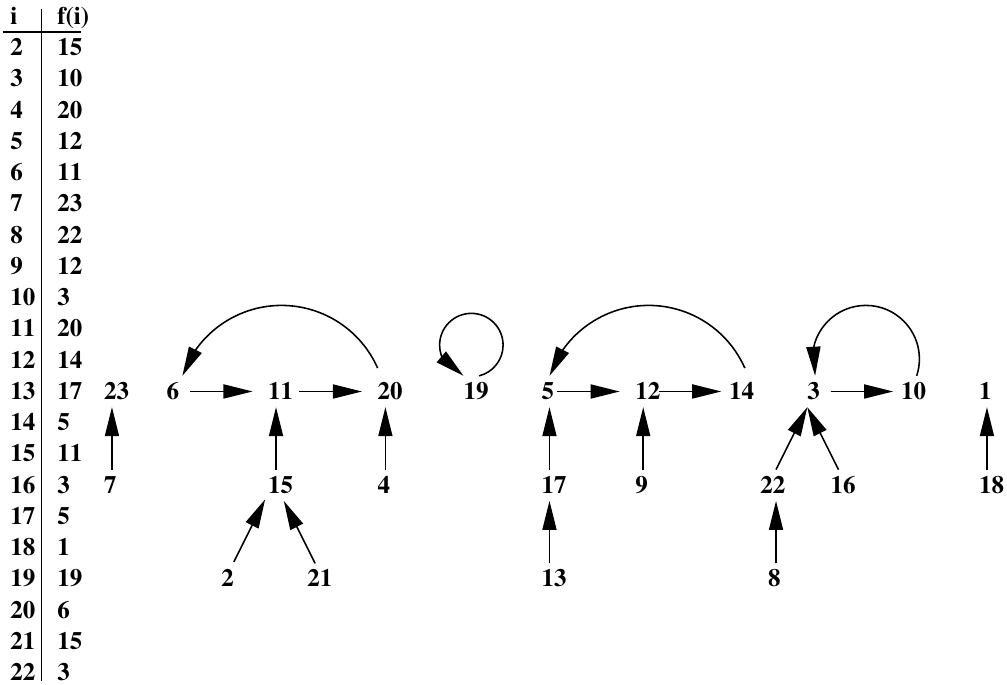}{The digraph of a function}

%\begin{figure}\label{digraph}
%\centerline{\hbox {\psfig {figure=Theta.ps,height=5in,width=5in}}}
%\caption{The digraph of a function}
%\end{figure}

This given, suppose that the digraph of $f$ is drawn as described above and
the cycles of $f$ are $c_1(f), \ldots, c_a(f)$, reading from left to right.
We let $r_{c_i(f)}$ and $l_{c_i(f)}$ denote the right and left endpoints of
the cycle $c_i(f)$ for $i= 1, \ldots,a$. Note that if $c_i(f)$ is a
1-cycle, then
we let $r_{c_i(f)} = l_{c_i(f)}$ be the element in the 1-cycle.
$\Theta(f)$ is obtained from $f$ by simply deleting the back edges
$(r_{c_i(f)},l_{c_i(f)})$ for $i = 1, \ldots, a$ and adding the directed
edges $(r_{c_i(f)},l_{c_{i+1}(f)})$ for $i =
1, \ldots, a-1$ plus the directed edges $(n,l_{c_1(f)})$ and
$(r_{c_a(f)},1)$.
That is, we remove all the back
edges that are above the line, and then we connect $n$ to the lefthand
endpoint of the first cycle, the righthand
endpoint of each cycle to the lefthand endpoint of the cycle following it,
and we connect the righthand endpoint of
the last cycle to $1$.
For example, $\Theta(f)$ is pictured in Figure~2  for the $f$ given in
Figure~1.
If there are no cycles in $f$, then $\Theta(f)$ is simply the result of
adding the directed edge $(n,1)$ to the digraph of $f$.

\fig{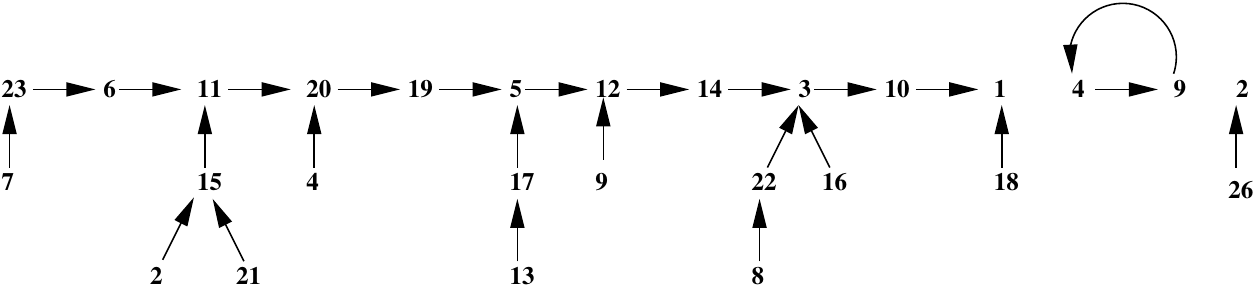}{$\Theta(f)$}

%\begin{figure}\label{image}
%\centerline{\hbox {\psfig {figure=Theta1.ps,height=4in,width=5in}}}
%\caption{$\Theta(f)$}
%\end{figure}

To see that $\Theta$ is a bijection, we shall describe how to define
$\Theta^{-1}$. The key observation is that we need only recover 
that the directed edges $(r_{c_i(f)},l_{c_{i+1}(f)})$ for $i =
1, \ldots, a-1$.  However it is easy to see that 
$r_{c_1(f)} = m_1$ is the largest element on the path from $n$ to $1$ 
in the tree $\Theta(f)$.  That is, $m_1$ is then largest element in its cycle and by definition, it is larger than all the largest elements in any other cycle so 
that $m_1$ must be the largest interior element on the path from $n$ to 1. 
Then by the same reasoning, $r_{c_2(f)} = m_2$ is the largest element 
on the path from $m_1$ to 1, etc. Thus we can find $m_1, \ldots, m_t$.    
More formally, given a tree  $T\in \vec{C}_{n,1}$, consider the path
$$m_0=n, x_1, \ldots, m_1, x_2, \ldots m_2, \ldots, x_t, \ldots , m_t, 1$$
where  $m_i$ is the maximum interior vertex on the path
from $m_{i-1}$ to $1$, $1\leq i\leq t$.  If $(m_{i-1},m_i)$ is an edge on
this path, then it is
understood that $x_i, \ldots , m_i = m_i$ consists of just one vertex and
we define $x_i=m_i$.  Note that by definition $m_0 = n > m_1 > \ldots > m_t$.
We obtain the digraph $\Theta^{-1}(T)$ from $T$ via the following
procedure.\\
\ \\
{\bf Procedure for computing $\Theta^{-1}(T):$}

\medskip\noindent
(1) First we declare that any  edge $e$ of $T$
which is not an edge of the path
from $n$ to $j$ is
an edge of $\Theta^{-1}(T)$.\\
\ \\
(2) Next we remove all edges of the form
$(m_t, 1)$  or $(m_{i-1}, x_i)$ for $1\leq i\leq t$.\\
\ \\
Finally for each $i$ with $1\leq i\leq t$, we consider the subpath
$x_i, \ldots , m_i$.\\
\ \\
(3) If $m_i=x_i$, create a directed loop
$(m_i, m_i)$.\\
\ \\
(4) If $m_i\neq x_i$,
convert the
subpath
$x_i, \ldots , m_i$ into the \\
directed cycle $x_i, \ldots, m_i, x_i$.\\
\ \\

Next we consider two important properties of the bijection $\Theta$.
First $\Theta$ has an important weight preserving property. 
We claim that if $\Theta(f) = T$, then 
\begin{equation}\label{weightpres}
q_nt_1 W(f) = W(T).
\end{equation}
That is, by our conventions, any backedge 
$(r_{c_{i}(f)},l_{c_{i}(f)})$ are descent edges so that its weight is
$q_{r_{c_{i}(f)}}t_{l_{c_{i}(f)}}$.  Thus the total weight of the backedges is
\begin{equation}\label{fnw}
\prod_{i=1}^a q_{r_{c_{i}(f)}}t_{l_{c_{i}(f)}}.
\end{equation}
Our argument above shows that all the new edges that we add are also
descent edges so that the weight of the new edges is
\begin{equation}
q_nt_{l_{c_1(f)}} (\prod_{i=1}^{a-1} q_{r_{c_{i}(f)}}t_{l_{c_{i+1}(f)}})
q_{r_{c_a(f)}}t_1 = q_nt_1 \prod_{i=1}^a q_{r_{c_{i}(f)}}t_{l_{c_{i}(f)}}.
\end{equation}
Since all the remaining edges have the same weight in both the digraph of
$f$ and in the digraph $\Theta(f)$, it follows that $q_nt_1W(f) = W(\Theta(f))$ as claimed.

It is easy to see that 
\begin{equation}
\sum_{f \in {\cal F}_n} W(f) = \prod_{i=2}^{n-1} [q_i(t_1 + \cdots + t_i) + 
p_i(s_{i+1} + \cdots + s_n)].
\end{equation}
Thus we have the following result which is implicit in \cite{ER1} and it explicit in \cite{RW}.
\begin{theorem}\label{weightthm1}
\begin{equation}
\sum_{T \in \vec{C}_{n,1}} W(T) = q_n t_1 \prod_{i=2}^{n-1} 
[q_i(t_1 + \cdots + t_i) + p_i(s_{i+1} + \cdots + s_n)].
\end{equation}
\end{theorem}

Next we turn to a second key property of the $\Theta$ bijection. It 
is easy to see from Figures 1 and 2 that 
deleting the back edges
$(r_{c_i(f)},l_{c_i(f)})$ for $i = 1, \ldots, a$ in $graph(f)$ and adding the directed edges $(r_{c_i(f)},l_{c_{i+1}(f)})$ for $i =
1, \ldots, a-1$ plus the directed edges $(n,l_{c_1(f)})$ and
$(r_{c_a(f)},1)$ to get $\Theta(f)$ does not change the indegree of 
any vertex except vertex 1. 
That is,    
\begin{equation}\label{Deg1}
indeg_{graph(f)}(i) = indeg_{\Theta(f)}(i)\ \mbox{for} \ i =2, \ldots ,n.
\end{equation}
It is also easy to see that in going from $graph(f)$ to $\Theta(f)$, the 
indegree of vertex 1 increases by 1, i.e., 
\begin{equation}\label{Deg2}
1+ indeg_{graph(f)}(i) = indeg_{\Theta(f)}(1).
\end{equation}
When we consider $\Theta(f)$ as an undirected graph $T$, then it is 
easy to see that $deg_T(i) = outdeg_{\Theta(f)} + indeg_{\Theta(f)}$. Thus 
since the outdegree of $i$ in $\Theta(f)$ is 1 if $ i \neq 1$ and the outdegree of 1 in $\Theta(f)$ is zero, equations (\ref{Deg1}) and (\ref{Deg2}) imply 
the following Theorem.

\begin{theorem}\label{degthm}
 Suppose that $T$ is the undirected tree corresponding 
to $\Theta(f)$ where $f \in {\cal F}_n$, then for $i =1, \ldots, n$, 
\begin{equation}\label{Deg3}
deg_T(i) = 1 + |f^{-1}(i)|.
\end{equation}
\end{theorem}
{\em Proof} By our definition of $graph(f)$, it follows that 
$indeg_{graph(f)}(i) = |f^{-1}(i)|$ for $i =1, \ldots, n$. 
Thus by (\ref{Deg1}), for $i = 2, \ldots, n$,
\begin{eqnarray*}
deg_T(i) &=&  outdeg_{\Theta(f)}(i) + indeg_{\Theta(f)}(i) \\
&=& 1 + indeg_{\Theta(f)}(i) \\
&=& 1 + indeg_{graph(f)}(i) \\
&=& 1 + |f^{-1}(i)|.
\end{eqnarray*}
Similarly by (\ref{Deg2}), 
\begin{eqnarray*}
deg_T(1) &=&  outdeg_{\Theta(f)}(1) + indeg_{\Theta(f)}(1) \\
&=& 0 + indeg_{\Theta(f)}(1) \\
&=& 1 + indeg_{graph(f)}(1) \\
&=& 1 + |f^{-1}(1)|.
\end{eqnarray*}
\hfill $\Box$

\section{Construction of the spanning forests from the the function
table in time $O(n)$}

In this section, we shall briefly outline the proof that one can 
compute the bijections  $\Theta$ and its inverse in linear time.  
Suppose we are given $f \in {\cal F}_n$. 
Our basic data structure for the  function $f$ is a list of pairs $
\langle i,f(i) \rangle$ for $i = 2, \ldots, n-1$.   
Our goal is to construct the directed 
graph of $\Theta(f)$ from our data structure for $f$, that is, 
for $i = 1, \ldots, n$, we want to find the set of pairs, $\langle i,t_i \rangle$,  such that there is directed edge from $i$ to $t_i$ in $\Theta(f)$. We shall prove the following. 

\begin{theorem}\label{linbij}
We can compute the bijection
$\Theta:{\cal F}_n \rightarrow \vec{C}_{n,1}$ 
and its inverse in linear time.
\end{theorem}
{\em Proof.}
We shall not try to give the most efficient algorithm to construct 
$\Theta(f)$ from $f$. Instead, we shall give an 
outline the basic procedure which shows that one can 
construct $\Theta(f)$ from $f$ in linear time. For 
ease of presentation, we shall organize our procedure so that it 
makes four linear time passes through the basic data structure for $f$ to 
produce the data structure for $\Theta(f)$.\\ 
\ \\
{\bf Pass 1.} {\it Goal: Find, in linear time in $n$, a set of representatives
$t_1, \ldots, t_r$ of the cycles of the directed graph of the function $f$.} \\
To help us find $t_1, \ldots, t_r$, we shall maintain an array 
$A[2], A[3], \cdots A[n-1]$,  
where for each $i$, $A[i]=(c_i, p_i, q_i)$ is a triple of integers such 
$c_i \in \{0,\ldots, n-1\}$ and $\{p_i,q_i\} \subseteq 
\{-1,2,\cdots, n-1\}$.  
The $c_i$'s will help us keep track of what loop we are in relative to 
the sequence of operations described below. Then our idea is to  
maintain, through the $p_i$ and $q_i$, a doubly linked list of the
locations $i$ in $A$ where $c_i=0$, 
and we obtain pointers to the first and last elements of this
doubly linked list. It is a standard exercise that these 
data structures can be maintained in
linear time.

Initially, all the $c_i$'s will be 
zero.  In general, if $c_i=0$, then $p_i$ will be the largest integer $j$ such that  
$2 \leq j < i$ for which $c_{j}=0$ if there is such a $j$ and 
$p_i = -1$ otherwise. 
Similarly, we set 
$q_i>i$ to be the smallest integer $k$ such that $n-1 \geq k > i$  for which
$c_{k}=0$ if there is such a $k$ and 
$q_i=-1$ if there is no such $k$.   If $c_{2}>0$, then $q_{2}$ is
the smallest integer $j > 2$  such that $c_{j}=0$ and 
${q_{2}}=-1$ if there is  no such integer $j$.  If $c_{n-1}>0$,  then
${p_{n-1}}$ is the largest integer $k < n_1$ such that $c_{k}=0$ and
${p_{n-1}}=-1$ if there is  no such integer $k$.

We initialize $A$ by setting $A[2]=(0,-1,q_{2})$, 
$A[i]=(0,i-1,i+1)$ for $m+1<i<n-1$, and $A[n-1]=(0,p_{n-1},-1)$.
If $2<n-1$ then $q_{2}=3$ and $p_{n-1}=n-2$.  Otherwise ($2=n-1$),
these quantities are both $-1$.\\
\ \\
{\bf LOOP(1):} Start with $i_1=2$, setting $c_{2}=1$. 
Compute   $f^0(2), f^1(2),  f^2(2), \ldots, f^{k_1}(2)$, each
time updating
$A$ by setting $c_{f^{j}(2)}=1$ and adjusting pointers, until, 
prior to setting $c_{f^{k_1}(2)}=1$, we discover that  either \\
\ \\
(1) $f^{k_1}(2) \in \{1, n\}$, in which case  
we have reached a node in ${graph}(f)$ which is not in the domain of $f$ and 
we start over again with the
$2$ replaced by the smallest $i$ for which  $c_i=0$, or \\
\ \\
(2) $x=f^{k_1}(2)$ already satisfies $c_x=1$.  This condition indicates that
the value $x$ has
already occurred in the sequence $2, f(2),  f^2(2), \ldots,
f^{k_1}(2)$.  Then we set $t_1=f^{k_1}(2)$.\\
\ \\
{\bf LOOP(2):} Start with $i_2=q_{m+1}$ which is the location of the first $i$ such that $c_i=0$, and repeat the calculation of LOOP1 with
$i_2$ instead of $i_1=2$.
In this manner, generate $f^0(i_2), f^1(i_2),  f^2(i_2), \ldots,
f^{k_2}(i_2)$, each time updating
$A$ by setting $c_{f^{j}(i_2)}=2$ and adjusting pointers,
until either \\
\ \\
(1) $f^{k_1}(i_2) \in \{1,n\}$, in which case we have 
reached a node in
${graph}(f)$ which is not in the domain of $f$ and we start
over again with the
$i_2$ replaced by the smallest $i$ for which  $c_i=0$, or \\
\ \\
(2) $x=f^{k_1}(i_2)$ already satisfies $c_x=2$. (This condition indicates that
the value $x$ has
already occurred in the sequence $i_2, f(i_2),  f^2(i_2), \ldots,
f^{k_1}(i_2)$.) Then we set $t_2=f^{k_1}(i_2)$.\\
\ \\
We continue this process until $q_{2}=-1$. 
At this point, we will have generated
$t_1, \ldots, t_r$, where the last loop was LOOP($r$). 
The array $A$ will be such that, for all $2\leq i\leq n-1$, 
$1\leq c_i\leq r$ identifies the LOOP in which
that particular domain value $i$ occurred in our computation described above.\\
\ \\
{\bf Pass 2.} {\em Goal: For $i =1, \ldots, r$, find the largest element $m_i$ in the cycle determined by $t_i$.}\\ 
It is easy to see that this computation can be  done in linear time by
one pass through the array $A$ computed in Pass 1 above. At the end of 
Pass 2, we set $l_i = f(m_i)$. Thus when we draw the cycle containing $t_i$ according to our definition of $\Theta_j(f)$, $m_i$ will be right most element in the and $l_i$ will be the left most element of the cycle containing $t_i$. However, at this point, we have not ordered the cycles appropriately.  This ordering will be done in the next pass. \\
\ \\
{\bf Pass 3.} {\em Goal:  Sort $(l_1,m_1), \ldots, (l_k,m_k)$ so that they are
appropriately ordered according the criterion for the bijection $\Theta(f)$ 
as described in by condition (a) -(d)}.\\
\ \\
Since we order the cycles from left to right according to decreasing 
maximal elements, it is then easy to see that our desired 
ordering can be constructed via a lexicographic bucket sort. 
(See Williamson's book \cite{W} for details on the fact that a 
lexicographic bucket sort can be carried out in linear time.)\\
\ \\ 
{\bf Pass 4.} {\em Goal: Construct the digraph of $\Theta(f)$ from 
the digraph of $f$.}\\
\ \\
We modify the table for $f$ to produce the table for $\Theta(f)$ as follows.
Assume that $(l_1,m_1), \ldots, (l_k,m_k)$ is the sorted list
coming out of Pass 3.
Then we modify the table for $f$ so that we add entries for 
the directed edges $\langle n, l_1\rangle $ and $ \langle m_k, 1 \rangle$ and
modify entries of the pairs starting with $m_1, \ldots, m_k$ so that their corresponding second elements are $l_2, \ldots, l_k, 1$ respectively.
This can be done in linear time using our data structures.

Next, consider the problem of computing the inverse of $\Theta$.
Suppose that we are given the data structure of the tree  $T \in \vec{C}_1$, i.e. we are given a 
set of pairs , $\langle i,t_i\rangle$,  such that there is a directed edge from $i$ to $t_i$ in $T$.
Recall that the computation of $\Theta^{-1}(T)$ consists of two basic steps. \\
\ \\
{\bf Step 1.} Given a tree  $T \in \vec{C}_{n,1}$, consider the path
$$m_0=n, x_1, \ldots, m_1, x_2, \ldots m_2, \ldots, x_t, \ldots , m_t, 1$$
where  $m_i$ is the maximum interior vertex on the path
from $m_{i-1}$ to $1$, $1\leq i\leq t$.  If $(m_{i-1},m_i)$ is an edge on
this path, then it is
understood that $x_i, \ldots , m_i = m_i$ consists of just one vertex and
we define $x_i=m_i$.  Note that by definition $m_0 = n > m_1 > \ldots > m_t$.\\
\ \\
First it is easy to see that by making one pass through the data structure for $F$, 
we can construct the directed path $n \rightarrow a_1 \rightarrow \ldots 
\rightarrow a_r$ where $1 = a_r$. 
In fact, we can construct a doubly linked list $(n,a_1, \ldots, a_{r-1},1)$ with pointers 
to the first and last elements in linear time. If we traverse the list in reverse order, 
$(1, a_{r-1}, \ldots, a_1,n)$, then it easy to see that $m_t = a_{r-1}$, $m_{t_1}$ is the 
next element in the list $(a_{r-2}, \ldots ,a_1)$ which is greater than $m_t$ and, in 
general, having found $m_i = a_s$, then $m_{i-1}$ is the first element in the list 
$(a_{s-1}, \ldots, a_1)$ which is greater than $m_i$.  Thus it is not difficult to see 
that we can use our doubly linked list to produce the factorization 
$$m_0=n, x_1, \ldots, m_1, x_2, \ldots m_2, \ldots, x_t, \ldots , m_t, 1 $$  
in linear time. \\
\ \\ 
{\bf Step 2.}
We obtain the digraph $\Theta^{-1}(T)$ from $T$ via the following
procedure.\\
\ \\
{\bf Procedure for computing $\Theta_j^{-1}(F):$}

\medskip\noindent
(1) First we declare that any  edge $e$ of $T$
which is not an edge of the path
from $n$ to $1$ is
an edge of $\Theta^{-1}(T)$.\\
\ \\
(2) Next we remove all edges of the form
$(m_t, 1)$  or $(m_{i-1}, x_i)$ for $1\leq i\leq t$.\\
\ \\
Finally for each $i$ with $1\leq i\leq t$, we consider the subpath
$x_i, \ldots , m_i$.\\
\ \\
(3) If $m_i=x_i$, create a directed loop
$(m_i, m_i)$.\\
\ \\
(4) If $m_i\neq x_i$, then ,
convert the
subpath
$x_i, \ldots , m_i$ into the \\
directed cycle $x_i, \ldots, m_i, x_i$.\\
\ \\
Again it is easy to see that we can use the data structure for 
$T$, our doubly linked list, and our path factorization, 
$m_0=n, x_1, \ldots, m_1, x_2, \ldots m_2, \ldots, x_t, \ldots , m_t, j$ to construct the 
data structure for $graph(f)$  where $f = \Theta^{-1}(T)$ in linear time.  
\hfil $\Box$.

Given that we can carry out the bijection $\Theta$ and its inverses in 
linear time, it follows that in linear time, we can reduce the problem of constructing  ranking and unranking algorithms for $\vec{C}_{n,1}$ to the problem 
of constructing ranking and unranking algorithms for the corresponding function class ${\cal F}_n$.  
In the next section, we will construct our desired  
ranking and unranking 
algorithms for the function classes corresponding to sets of trees 
$\vec{C}_{n,\vec{s}}$ and $\vec{C}_{n,S}$ described in the introduction.

\section{Ranking and Unranking Algorithms for 
Trees with a Fixed Degree Sequence.}

Recall that if $T \in C_n$, then 
$\sum_{i=1}^n deg_T(i) = 2n-2$. Let 
$\vec{s} = \langle s_1, \ldots, s_n \rangle$ be any sequence of positive integers such that $\sum_{i=1}^n s_i =2n-2$. Then we define $\vec{C}_{n,\vec{s}} = \{T \in \vec{C}_{n,1}: \langle deg_T(1), \ldots, deg_T(n)\rangle = \vec{s}\}$. Similarly if $S =\{1^{\alpha_1}, \ldots, (n-1)^{\alpha_{n-1}}\}$ is a multiset such that $\sum_{i=1}^{n-1} \alpha_i \cdot i = 2n-2$ and $\sum_{i=1}^n \alpha_i = n$, then we define   
$\vec{C}_{n,S} = \{T \in \vec{C}_{n,1}: 
\{deg_T(1), \ldots, deg_T(n)\} = S \}$. 
The main goal of this section is to construct $n^2log(n)$ time algorithms 
for ranking and unranking trees in $\vec{C}_{n,\vec{s}}$ and $\vec{C}_{n,S}$. 

So assume that $\vec{s} = \langle s_1, \ldots, s_n \rangle $ is a sequence of 
positive integers such that  
$\sum_{i=1}^n s_i = 2n-2$. By Theorem \ref{degthm}, it follows that if 
$\Theta(f)  = T$, then $deg_T(i) = 1 +|f^{-1}(i)|$ for $i =1, \ldots n$. 
It follows that 
\begin{equation}\label{vecs}
\Theta^{-1}(\vec{C}_{n, \vec{s}}) = 
\{f \in {\cal F}_n: \langle |f^{-1}(1)|, \ldots, |f^{-1}(n)|\rangle 
= \langle s_1 -1, \ldots, s_n-1\rangle\}.
\end{equation}
Since a function $f \in {\cal F}_n$ is clearly determined by 
the sequence $\langle f^{-1}(1), \ldots, f^{-1}(n)\rangle$, it follows 
from our results in Section 2 that the problem of finding an algorithms to rank and unrank trees $\vec{C}_{n,\vec{s}}$ can be reduced to the problem 
of ranking and unranking ordered set partitions 
$\langle \pi_1, \ldots, \pi_n \rangle$ of $\{2, \ldots, n-1\}$ where the sizes of the sets are specified.  
That is, we need to find an algorithm to 
rank and unrank ordered set partitions in $\Pi_{n,\vec{s}}$, the set of all sequences of pairwise disjoint sets 
$\langle \pi_1, \ldots, \pi_n \rangle$ such that 
$\bigcup_{k=1}^n \pi_k= \{2, \ldots, n-1\}$ and 
$|\pi_k| = s_k -1$ for $k =1, \ldots, n$. The total number of elements in 
$\Pi_{n,\vec{s}}$ is clearly the multinomial coefficient 
$\binom{n-2}{s_1-1, \ldots, s_n -1} = \binom{n-2}{s_1-1} \binom{n-2 -(s_1-1)}{s_2 -1} \cdots \binom{n-2-(\sum_{i=1}^n s_i-1)}{s_n-1}$. Thus our 
first step is to develop a simple algorithm to rank and unrank objects corresponding to a product of binomial coefficients 
$\prod_{i=1}^k \binom{a_i}{b_i}$.  

For a single binomial coefficient $\binom{n}{k}$, we shall rank and unrank the set ${\cal DF}_{n,k}$ of decreasing functions $f:\{1, \ldots, k\} \rightarrow \{1, \ldots, n\}$ relative to lexicographic order. A number of authors have developed ranking and unranking algorithms for ${\cal DF}_{n,k}$.  We shall follow the method of Williamson \cite{W}. First, we identify a function 
$f:\{1, \ldots, k\} \rightarrow \{1, \ldots, n\}$ with the decreasing sequence 
$\langle f(1), \ldots, f(k) \rangle$ where 
$n \geq f(1) > \ldots > f(k) \geq 1$.  We can then think of the sequences as specifying a node in a planar tree $T_{{\cal DF}_{n,k}}$ which can be constructed recursively as follows. At level 1, 
the nodes of $T_{{\cal DF}_{n,k}}$ are labeled   
$k, \ldots, n$ from left to right specifying the choices for $f(1)$. Next below a node $j$ at level one, we attach a tree corresponding to 
$T_{{\cal DF}_{j-1,k-1}}$ where a tree  $T_{{\cal DF}_{a,1}}$ consists of a tree with a single vertex labeled $a$. Figure 3 pictures the tree 
$T_{{\cal DF}_{6,3}}$.

\fig{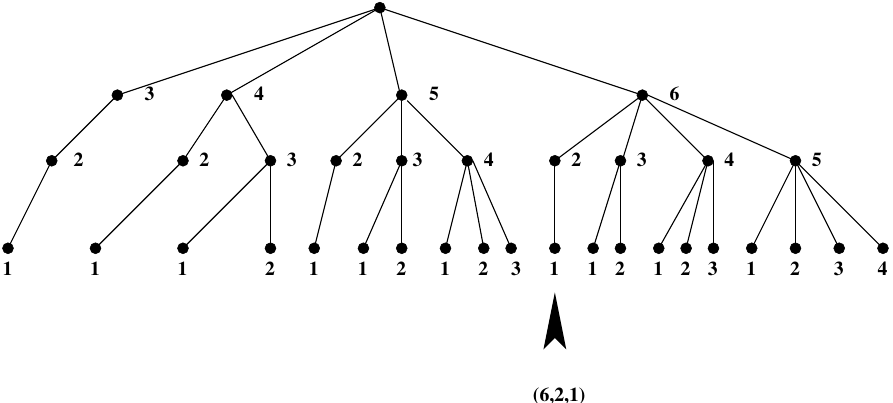}{The tree $T_{{\cal DF}_{6,3}}$}

%\begin{figure}[htbp]
%\vspace{1.5in}
%\fcaption{Labeled tree {\it T.}}
%\end{figure}
%\begin{center}
%\
%\psfig{figure=T73.eps,width=3in,height=1.5in}
%\fcaption{The tree $T_{{\cal DF}_{6,3}}$}
%\end{center}
%\end{figure}

Then the decreasing sequence (6,2,1) corresponds to the node which is specified with an arrow.  It is clear that the sequences corresponding to the 
nodes at the bottom of the tree $T_{{\cal DF}_{6,3}}$ appear in lexicographic order from left to right.  Thus the rank of any sequence $\langle 
f(1), \ldots, f(k) \rangle \in {\cal DF}_{n,k}$ is the number of nodes at the bottom of the tree to the left of the node corresponding to $\langle 
f(1), \ldots, f(k) \rangle$. Hence the sequence $\langle 6,2,1\rangle$
 has rank 10 in the tree $T_{{\cal DF}_{6,3}}$. 

This given, suppose we are given a sequence $\langle f(1), \ldots, f(k) \rangle$ in $T_{{\cal DF}_{n,k}}$. Then the number of leaves in the subtrees 
corresponding the nodes $k, \ldots, f(1) -1$ are respectively 
$\binom{k-1}{k-1}, \binom{k}{k-1}, \ldots, \binom{f(1)-2}{k-1}$. Thus the total number of leaves in those subtrees is 
$$\binom{k-1}{k-1} + \binom{k}{k-1}+ \cdots + \binom{f(1)-2}{k-1} = 
\binom{f(1)-1}{k}.$$
Here we have used the well known identity that 
$\sum_{s=k-1}^{t-1} \binom{s}{k-1} = \binom{t}{k}$. 
It follows that the rank of 
$\langle f(1), \ldots, f(k) \rangle$ in $T_{{\cal DF}_{n,k}}$ equals 
$\binom{f(1)-1}{k}$ plus the rank of $\langle f(2), \ldots, f(k)\rangle$ in 
$T_{{\cal DF}_{f(1)-1,k-1}}$.  The following result, stated in 
\cite{W}, then easily follows by induction.

\begin{theorem}\label{rankdecfn}
Let $f:\{1, \ldots, k\} \rightarrow \{1, \ldots, n\}$ be a descreasing function. Then the rank of $f$ relative to the lexicographic order on ${\cal DF}_{n,k}$ is 
\begin{equation}\label{decfnrank}
rank_{{\cal DF}_{n,k}}(f) = \binom{f(1)-1}{k} + \binom{f(2)-1}{k-1} + \cdots + \binom{f(k)-1}{1}.
\end{equation}
\end{theorem}

It is then easy to see from Theorem \ref{rankdecfn}, that the following 
procedure, as described by Williamson in \cite{W}, gives the unranking procedure for ${\cal DF}_{n,k}$.\\
\begin{theorem}\label{unrankdecfn}
The following procedure $UNRANK(m)$ computes 
$$f =\langle f(1),\ldots, f(k)\rangle$$
 such that $Rank_{{\cal DF}_{n,k}}(f) =m$ for any 
$1 \leq k \leq n$ and $0 \leq m \leq \binom{n}{k}-1$.\\
\ \\
{\bf Procedure} UNRANK(m) \\
\begin{description}
\item[initialize] $m':=m$, $t:=1$, $s:= k$; ($1 \leq k \leq n$, $0 \leq m \leq \binom{n}{k}-1$)
\item[while] $t \leq k$ \ {\bf do}
\begin{description}
\item[begin]
\begin{description}
\item $f(t)-1 = max\{y: \binom{y}{s} \leq m'\}$;
\item $m' := m' - \binom{f(t)-1}{s}$;
\item $t := t+1$;
\item $s := s-1$;
\end{description}
\item[end]
\end{description}
\end{description}

\end{theorem}

Note that 
$|{\cal DF}_{a_1,b_1} \times {\cal DF}_{a_2,b_2} \times \cdots \times {\cal DF}_{a_t,b_t}| = \prod_{i=1}^t \binom{a_i}{b_i}$.
Thus we can use ${\cal DF}_{a_1,b_1} \times {\cal DF}_{a_2,b_2} \times \cdots \times {\cal DF}_{a_t,b_t}$ as the set of objects corresponding to a product 
of binomial coefficients. We shall idenitfy an element 
$$(f_1, \ldots, f_t) \in  
{\cal DF}_{a_1,b_1} \times {\cal DF}_{a_2,b_2} \times \cdots \times {\cal DF}_{a_t,b_t}$$ with a sequence 
$$\langle f_1(1), \ldots, f_1(b_1), f_2(1), \ldots, f_2(b_2), \ldots, f_t(1), \ldots, f_t(b_t)\rangle$$ 
and rank these sequences according to lexicographic order. To define our ranking and unranking proceedure for this set of sequences, we first 
need to define a product relation on planar trees.

Given a rooted planar tree $T$, let $L(T)$ be the numbers of leaves of 
$T$ and $Path(T)$ be the set of paths which go from the root to a leaf.  
Then for any path $p \in Path(T)$, we define the rank of $p$ relative to $T$, 
$rank_T(p)$, to be the number of leaves of $T$ that lie to the left of $p$.

Given two rooted planar trees $T_1$ and $T_2$, $T_1 \otimes T_2$ is the 
tree that results from $T_1$ by replacing each leaf of $T_1$ 
by a copy of $T_2$, see Figure 3. If the vertices of $T_1$ and $T_2$ 
are labeled, then we shall label the vertices of $T_1 \otimes T_2$ according 
to the convention that each vertex $v$ in $T_1$ have the same label in 
$T_1 \otimes T_2$ that it has in $T_1$ and each vertex $w$ in a copy of 
$T_2$ that is decendent from a leaf labeled $l$ in $T_1$ has a label 
$(l,s)$ where $s$ is the label of $w$ in $T_2$. 
Given rooted planar trees $T_1, \ldots, T_k$ where $k \geq 3$, we can define 
$T_1 \otimes T_2 \otimes \cdots \otimes T_k$ by induction as 
$(T_1 \otimes \cdots \otimes T_{k-1}) \otimes T_k$. Similarly if 
$T_1 \ldots, T_k$ are labeled rooted planar trees, we can define 
the labeling of $T_1 \otimes T_2 \otimes \cdots \otimes T_k$ by the same 
inductive process. 

\fig{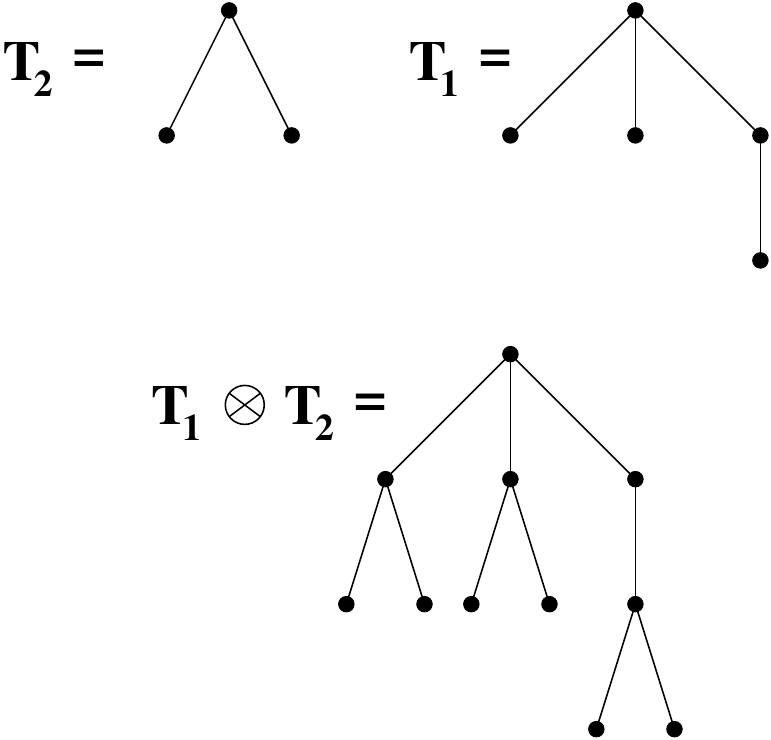}{The operation $T_1 \otimes T_2$.}

%\begin{figure}\label{otimes}
%\centerline{\hbox {\psfig {figure=prod.ps,height=3in,width=3in}}}
%\caption{The operation $T_1 \otimes T_2$.}
%\end{figure}

Now suppose that we are given two rooted planar trees $T_1$ and $T_2$ and 
suppose that $p_1 \in Path(T_1)$ and $p_2 \in Path(T_2)$.  Then we 
define the path $p_1 \otimes p_2$ in $T_1 \otimes T_2$ which follows 
$p_1$ to its leaf $l$ in $T_1$ and then follows $p_2$ in the copy of 
$T_2$ that sits below leaf $l$ to a leaf $(l,l')$ in $T_1 \otimes T_2$. 
Similarly, given paths $p_i \in T_i$ for $i = 1, \ldots k$, 
we can define a path  $p = p_1 \otimes \cdots \otimes p_k \in Path(T_1 \otimes T_2 \otimes \cdots \otimes T_k)$ by induction as 
$(p_1 \otimes \cdots \otimes p_{k-1}) \otimes p_k$.

Next we give two simple lemmas that tell us how to rank and unrank the set 
of paths in such trees.

\begin{lemma}\label{prodrank}  Suppose that $T_1, \ldots ,T_k$ are rooted planar trees and $T = T_1 \otimes T_2 \otimes \cdots \otimes T_k$.  Then for any path 
$p = p_1 \otimes \cdots \otimes p_k \in Path(T)$, 
\begin{equation}\label{prod-rank}
rank_T(p) = \sum_{j=1}^k rank_{T_j}(p_j) \prod_{l=j+1}^k L(T_{l})
\end{equation} 
\end{lemma}
{\em Proof}.  We proceed by induction on $k$. Let us assume that 
$T_1, \ldots T_k$ are labeled rooted planar trees. 

First suppose that $k =2$ and that $p_1$ is a path that goes from 
the root of $T_1$ to a leaf labeled $1_1$ and $p_2$ goes from the root 
of $T_2$ to a leaf labeled $l_2$.  Thus $p_1 \otimes p_2$ goes from 
the root of $T_1 \otimes T_2$ to the leaf $l_1$ in $T_1$ and then proceeds 
to the leaf $(l_1,l_2)$ in $T_1 \otimes T_2$. Now for each 
leaf $l'$ to the left of $l_1$ in $T_2$, there are $L(T_2)$ leaves 
of $T_1 \otimes T_2$ that lie to left of $(l_1,l_2)$ coming from 
the leaves of the copy of $T_2$ that sits below $l'$. Thus 
there are a total of $L(T_2) \cdot rank_{T_1}(p_1)$ such leaves. The only 
other leaves of $T_1 \otimes T_2$ that lie to left of $p_1 \otimes p_2$ 
are the leaves of the form $(l_1,l'')$ where $l''$ is to left of $p_2$ in 
$T_2$.  There are $rank_{T_2}(p_2)$ such leaves.  Thus there are a total of 
$rank_{T_2}(p_2) + L(T_2) \cdot rank_{T_1}(p_1)$ leaves to left of 
$p_1 \otimes p_2$ and hence 
$$rank_{T_1 \otimes T_2}(p_1 \otimes p_2) = 
rank_{T_2}(p_2) + L(T_2) \cdot rank_{T_1}(p_1) $$
as desired. 

Next assume that (\ref{prod-rank}) holds for $k < n$ and that $n \geq 3$. Then 
\begin{eqnarray*}
&&rank_{T_1 \otimes \cdots \otimes T_n}(p_1 \otimes \cdots \otimes p_n) = \\
&&rank_{(T_1 \otimes \cdots \otimes T_{n-1}) \otimes T_n}
(p_1 \otimes \cdots \otimes p_{n-1}) \otimes p_n) = \\
&& rank_{T_n}(p_n) + L(T_n) 
(\sum_{j=1}^{n-1} rank_{T_j}(p_j) \prod_{l=j+1}^{n-1} L(T_{l})) =\\
&&\sum_{j=1}^n rank_{T_j}(p_j) \prod_{l=j+1}^n L(T_{l}).
\end{eqnarray*}
\hfil$\Box$.

This given, it is easy to develop an algorithm for unranking in a product 
of trees.  The proof of this lemma can be found in \cite{W}. 

\begin{lemma}\label{produnrank}  
Suppose that $T_1, \ldots ,T_k$ are rooted planar trees and $T = T_1 \otimes T_2 \otimes \cdots \otimes T_k$.  Then given a  
$p \in Path(T)$ such that $rank_T(p) = r_0$, 
$p = p_1 \otimes \cdots \otimes p_k \in Path(T)$ where 
$rank_{T_i}(p_i) = q_i$ and 
\begin{eqnarray}
r_0 &=& q_1 \prod_{l=2}^k L(T_{l}) + r_1 \ 
\mbox{where $0 \leq r_1 <  \prod_{l=2}^k L(T_{l})$}, \\
r_1 &=& q_2 \prod_{l=3}^k L(T_{l}) + r_2 \ \mbox{where $0 \leq r_2 <  
\prod_{l=3}^k L(T_{l})$},\\
&\vdots& \\
r_{k-2} &=& q_{k-1}L(T_k) + r_{k-1} \ \mbox{where $0 \leq r_{k-1} < L(T_k)$} \ \mbox{and} \\
r_{k-1} &=& q_k 
\end{eqnarray} 
\end{lemma}

It follows that ranking and unranking our sequences  
$$\langle f_1(1), \ldots, f_1(b_1), f_2(1), \ldots, f_2(b_2), \ldots, f_t(1), \ldots, f_t(b_t)\rangle$$ 
coresponding to an element $(f_1, \ldots, f_t) \in {\cal DF}_{a_1,b_1} \times {\cal DF}_{a_2,b_2} \times \cdots \times {\cal DF}_{a_t,b_t}$, we need only rank and unrank the leaves with respect to the tree
\begin{equation}\label{bigtree}
T_{(a_1, b_1), \ldots, (a_t,b_t)} = T_{{\cal DF}_{a_1,b_1}} \otimes T_{{\cal DF}_{a_2,b_2}} \otimes \cdots \otimes T_{{\cal DF}_{a_t,b_t}}.
\end{equation}
That is, consider a path $p = p_1 \otimes p_2 \otimes \cdots \otimes p_t \in T_{(a_1, b_1), \ldots, (a_t,b_t)}$.  For each $i$,  
$p_i$ corresponds to a sequence $\langle f_i(1), \ldots, f_i(b_i)\rangle$ in  
${\cal DF}_{a_i,b_i}$ and hence $p$ corresponds to the sequence 
$$\langle f_1(1), \ldots, f_1(b_1), f_2(1), \ldots, f_2(b_2), \ldots, f_t(1), \ldots, f_t(b_t)\rangle.$$

We are now in position to give an algorithm to 
rank and unrank ordered set partitions in $\Pi_{n,\vec{s}}$, the set of all sequences of pairwise disjoint sets 
$\langle \pi_1, \ldots, \pi_n \rangle$ such that 
$\bigcup_{k=1}^n \pi_k= \{2, \ldots, n-1\}$ and 
$|\pi_k| = s_k -1$ for $k =1, \ldots, n$. Since the 
total number of elements in 
$\Pi_{n,\vec{s}}$ is clearly the multinomial coefficient 
$\binom{n-2}{s_1-1, \ldots, s_n -1} = \binom{n-2}{s_1-1} \binom{n-2 -(s_1-1)}{s_2 -1} \cdots \binom{n-2-(\sum_{i=1}^n s_i-1)}{s_n-1}$, we shall identify 
an ordered set partition $\pi = \langle \pi_1, \ldots, \pi_n \rangle$ with an element 
$$(f_1, \ldots, f_n) \in {\cal DF}_{n-2,s_1-1} \times 
{\cal DF}_{n-2 -(s_1-1),s_2 -1} \times \cdots \times 
{\cal DF}_{n-2-(\sum_{i=1}^n s_i-1),s_n-1}$$
as follows. 
Suppose that $n=12$ and 
$\vec{s} = (s_1,\ldots, s_{12}) =     (1,1,3,1,4,1,3,1,2,1,3,1)$. 
Note that $\sum_{i=1}^{12} s_i = 2(12)-2 = 22$ so that this is a possible 
degree sequence for a tree in $\vec{C}_{12}$. For example, the tree 
$T_0 \in \vec{C}_{12,1}$ pictured in Figure 5  has this degree sequence 
when considered as a tree.

\fig{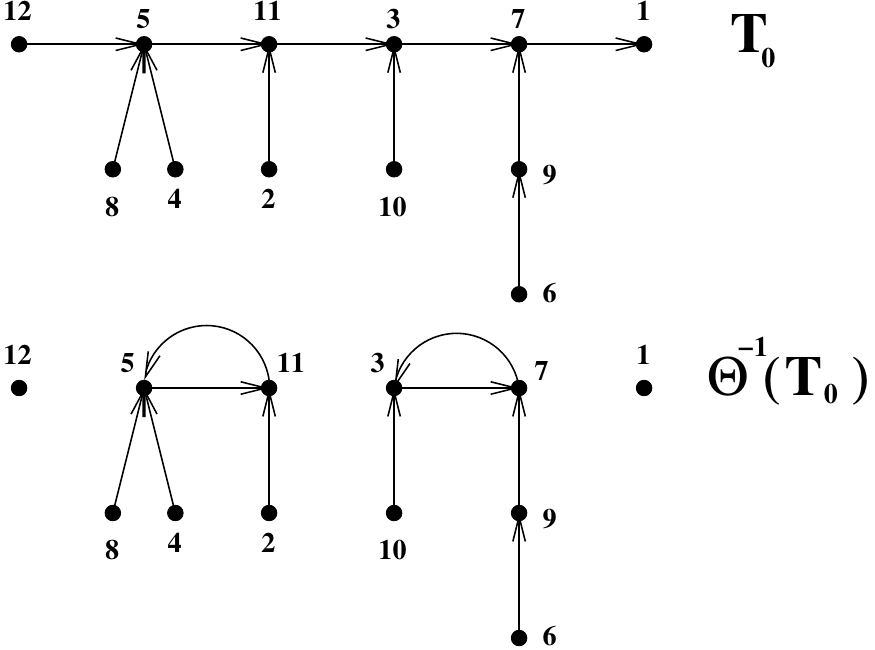}{The tree $T_0$}

%\begin{figure}[htbp]
%\vspace{1.5in}
%\fcaption{Labeled tree {\it T.}}
%\end{figure}
%\begin{center}
%\
%\psfig{figure=T0.eps,width=2in,height=1in}
%\fcaption{The tree $T_0$}
%\end{center}
%\end{figure}

Note that in this case,   
$\langle s_1-1,\ldots, s_{12}-1\rangle = \langle 0,0,2,0,3,0,2,0,1,0,2,0
\rangle$. 
Also in Figure 5, we have pictured the graph of 
$f= \Theta^{-1}(T_0)$ and in this 
case 
$$\pi_{T_0} = \langle f^{-1}(1), \ldots, f^{-1}(12)\rangle  = \langle \emptyset,\emptyset, \{10,7\}, \emptyset, \{11,8,4\}, \emptyset, 
\{9,3\}, \emptyset, \{6\}, \emptyset, \{5,2\}, \emptyset\rangle.$$
It will be more efficient for our ranking and unranking procedure to 
order the set partition by increasing size of the parts. Thus we will 
make one pass through the to extract the $|f^{-1}(i)|$ for each $i$ and 
its relative rank for those parts of the same size. In this case, 
we would produce the following list \\
{\tiny $\langle (1,0,0),(2,0,1),(3,2,0),(4,0,2),(5,3,0),(6,0,3),(7,2,1),(8,0,4),(9,1,0), (10,0,5),(11,2,2),(12,0,6)\rangle$.}\\
Here for example, then entry $(7,2,1)$ means that the size of $f^{-1}(7)$ is 
2 and there is one element $i < 7$ such that $|f^{-1}(i)| =2$.  We can 
then do a lexicographic bucket sort to produce a list of the elements according to lexicographic order on the last two entries of this list in linear time, 
see \cite{W}. Thus in our example we would produce the following list. \\
{\scriptsize $\langle (1,0,0),(2,0,1),(4,0,2),(6,0,3),(8,0,4),(10,0,5),(12,0,6), (9,1,0), 
(7,2,1), (3,2,0), (11,2,2),(5,3,0)\rangle$.}\\
The set partition corresponding to this order is 
$$\langle \emptyset,\emptyset,\emptyset,\emptyset,\emptyset,\emptyset,\emptyset,\{6\}, \{10,7\},\{9,3\},\{5,2\},\{11,8,4\}\rangle.$$
We can ignore the $\emptyset$'s and just consider the reduced partition 
$$\overline{\pi}_{T_0} = \langle \{6\}, \{10,7\},\{9,3\},\{5,2\},\{11,8,4\}\rangle.$$
 More generally, let $D=\{2, \ldots , n-1\}$. 
Let $\overline{\pi} = (B_1, \ldots, B_k)$ be an ordered
set partition of $D$ where each block $B_i$ is nonempty and ordered in
decreasing order, $B_i=\{b_{i,1}> \ldots > b_{i,t_i}\}$, that comes from some 
element $\pi$ in $\Pi_{12,\vec{s}}$ as described above.
Let $(r_{1,1}, \ldots, r_{1,t_1})$ be the ordered sequence of ranks of
the respective $(b_{1,1}, \ldots, b_{1,t_1})$ in $D$.
In general, let $(r_{i,1}, \ldots, r_{i,t_i})$ be the ranks of the
respective $(b_{i,1}, \ldots, b_{i,t_i})$ in
$D\setminus \cup_{j=1}^{i-1}B_j$.
For $\overline{\pi}_{T_0}$,\\
$(r_{1,1}) = (5)$ can be considered an element of 
${\cal DF}_{10,1}$ \\
$(r_{2,1},r_{2,2}) = (8,5)$ can be considered an element of 
${\cal DF}_{9,2}$ \\
$(r_{3,1},r_{3,2}) = (6,2)$ can be considered an element of 
${\cal DF}_{7,2}$ \\
$(r_{4,1},r_{4,2}) = (3,1)$ can be considered an element of 
${\cal DF}_{5,2}$ \\
$(r_{5,1},r_{5,2},r_{5,3}) = (3,2,1)$ can be considered an element of 
${\cal DF}_{3,3}$\\ 
Thus we can think of $\overline{\pi}_{T_0}$ as the sequence 
$\langle 6,8,5,6,2,3,1,3,2,1 \rangle$ coming from an element 
of 
${\cal DF}_{10,1}\times {\cal DF}_{9,2} \times {\cal DF}_{7,2} \times 
{\cal DF}_{5,2}\times {\cal DF}_{3,3}$ or as a leaf in the tree  
$T_{{\cal DF}_{10,1}}\otimes T_{{\cal DF}_{9,2}} \otimes T_{{\cal DF}_{7,2}} 
\otimes T_{{\cal DF}_{5,2}} \otimes T_{{\cal DF}_{3,3}}$. 
Note that the size of the trees needed for the product lemma, Lemma 
\ref{prodrank}, are \\
$|T_{{\cal DF}_{3,3}}| = \binom{3}{3} = 1$,\\
$|T_{{\cal DF}_{5,2}} \otimes T_{{\cal DF}_{3,3}}| = \binom{5}{2} \cdot \binom{3}{3} = 10$,\\
$|T_{{\cal DF}_{7,2}} 
\otimes T_{{\cal DF}_{5,2}} \otimes T_{{\cal DF}_{3,3}}| = 
\binom{7}{2} \cdot \binom{5}{2} \cdot \binom{3}{3} = 210$,\\
$|T_{{\cal DF}_{9,2}} \otimes T_{{\cal DF}_{7,2}} 
\otimes T_{{\cal DF}_{5,2}} \otimes T_{{\cal DF}_{3,3}}
T_{{\cal DF}_{7,2}} 
\otimes T_{{\cal DF}_{5,2}} \otimes T_{{\cal DF}_{3,3}}| = 
\binom{9}{2} \cdot \binom{7}{2} \cdot \binom{5}{2} \cdot \binom{3}{3} = 
7560$,\\
$|T_{{\cal DF}_{10,1}}\otimes T_{{\cal DF}_{9,2}} \otimes T_{{\cal DF}_{7,2}} 
\otimes T_{{\cal DF}_{5,2}} \otimes T_{{\cal DF}_{3,3}}| =
\binom{10}{1} \binom{9}{2} \cdot \binom{7}{2} \cdot \binom{5}{2} \cdot \binom{3}{3} = 
75600$.\\
Thus we can apply Lemma \ref{prodrank} and conclude that the
\begin{eqnarray*}
rank(\overline{\pi}_{T_0}) &=& 
rank_{{\cal DF}_{10,1}}(\langle 5 \rangle) \times 7560 \\
&&+ rank_{{\cal DF}_{9,2}}(\langle 8,5 \rangle) \times 210\\
&&+ rank_{{\cal DF}_{7,2}}(\langle 6,2 \rangle) \times 10\\
&&+ rank_{{\cal DF}_{5,2}}(\langle 3,1 \rangle) \times 1.
\end{eqnarray*}
By Theorem \ref{rankdecfn}, we have that 
\begin{eqnarray*}
rank_{{\cal DF}_{10,1}}(\langle 5 \rangle) &=& \binom{5-1}{1} = 4,\\
 rank_{{\cal DF}_{9,2}}(\langle 8,5 \rangle) &=& \binom{8-1}{2} + 
\binom{5-1}{1} = 21 +4 = 25,\\
rank_{{\cal DF}_{7,2}}(\langle 6,2 \rangle)&=& \binom{6-1}{2} + \binom{2-1}{1} = 10 +1 = 11, \ \mbox{and}\\
rank_{{\cal DF}_{5,2}}(\langle 3,1 \rangle)&=& \binom{3-1}{2} + \binom{1-1}{1} = 1 +0 = 1.
\end{eqnarray*}
Thus $$rank(\overline{\pi}_0) = (4\times 7560) +
(25 \times 210)+ (11\times 10) + 
(1 \times 1) = 35,601.$$
Hence the tree $T_0$ pictured in Figure 5  has rank 35,601 among all 
the trees \\
$T \in \vec{C}_{12,\langle 0,0,2,0,3,0,2,0,1,0,2,0
\rangle}$.

If we are given, the degree sequence $\vec{s}$, we can assume that we preprocess the sizes of the trees needed to apply the product lemma, Lemma \ref{prodrank}.  Thus we need only compute $O(n)$ products, additions, 
and multinomial coefficients.  Again given, $\vec{s}$, we can construct a table of all the possible binomial coefficients that we need as part of the preprocessing. Thus to find 
the rank of tree $T_0$ requires only a linear number of muliplications, additions and table look ups for numbers 
$x < |\vec{C}_{n,\vec{s}}|$. Since each $x$ requires at most $O(nlog(n))$ 
bits, 
it is easy to see that these operations require at most $O(n^2log(n))$ bit 
operations. 
The only other contribution to the complexity of the algorithm is the 
time it takes go from the representation of the tree to the 
corresponding rank sequences  
$$\langle r_{1,1}, \ldots, r_{1,t_1}, r_{2,1}, \ldots, r_{2,t_2}, 
\ldots, r_{k,1}, \ldots, r_{k,t_k}\rangle.$$ It is easy to see from the fact that we can compute $\Theta^{-1}$ in linear time that we can start with a tree 
and produce the ordered set partition $\overline{\pi}_T = (B_1, \ldots, B_k)$ in linear time.  Thus to complete our analysis of the 
the complexity of the ranking prodeedure, we 
need to know the complexity of the transformation between the 
$(B_1, \ldots, B_k)$ and ranks 
$$\langle r_{1,1}, \ldots, r_{1,t_1},r_{2,1}, \ldots, r_{2,t_2}, 
\ldots, r_{k,1}, \ldots, r_{k,t_k}\rangle.$$
\begin{lemma}\label{rankthm}
Let $D=\{2, \ldots , n-1\}$. Let $(B_1, \ldots, B_k)$ be an ordered
partition of $D$ where each block $B_i$ is nonempty and ordered in
decreasing order, $B_i=(b_{i,1}, \ldots, b_{i,t_i})$.
Let $(r_{1,1}, \ldots, r_{1,t_1})$ be the ordered sequence of ranks of
the respective $(b_{1,1}, \ldots, b_{1,t_1})$ in $D$.
In general, let $(r_{i,1}, \ldots, r_{i,t_i})$ be the ranks of the
respective $(b_{i,1}, \ldots, b_{i,t_i})$ in
$D\setminus \cup_{j=1}^{i-1}B_j$.  Given the sequences
$(r_{i,1}, \ldots, r_{i,t_i})$, $i=1, \dots, k$, the sets partition 
$B_1, \ldots, B_k$ can be constructed in worst case time
$O(n^2log(n))$. Conversly, given the set partition $B_1, \ldots, B_k$
the sequences $(r_{i,1}, \ldots, r_{i,t_i})$, $i=1, \dots, k$ can be constructed in worst case time $O(n^2log(n))$.
\end{lemma}
{\it Proof}. First we can make one pass through the list and set 
$b_{i,j}: b_{i,j} -1$ so $(B_1, \ldots, B_k)$ become an ordered
partition of $\{1, \ldots, n-2\}$. Conversly, we can go from 
an ordered set partition $(B_1, \ldots, B_k)$ of $\{1, \ldots, n-2\}$ to an 
ordered set partition of $\{2, \ldots, n-1\}$ by setting 
$b_{i,j}: = b_{i,j}+1$. Thus there is no loss in assuming that 
$(B_1, \ldots, B_k)$ is an ordered
partition of $\{1, \ldots, n-2\}$.

This given, it will then be the case that 
$(b_{1,1}, \ldots, b_{1,t_1}) = (r_{i,1}, \ldots, r_{i,t_i})$. 
Then it will take $O(n)$ comparisions of numbers less than or equal to 
$n$ to construct 
a sequence $f(1), \ldots, f(n)$ where $f(i) =|\{b_1,j: b_1,j  < i,j=1, \ldots, t_1\}|$.
Then it will take $nlog(n)$ steps to create the sequences 
$(\bar{b}_{i,1}, \ldots, \bar{b}_{i,t_i})$ for $i =2, \ldots,k$ 
where $\bar{b}_{i,j} = b_{i,j} - f(b_{i,j})$.  
It then easily follows 
that we have reduced the problem to finding the transformation from the ranks 
$(r_{i,1}, \ldots, r_{i,t_i})$, $i=2, \dots, k$, to set partitions 
$\bar{B}_2, \ldots, \bar{B}_k$ which we can do by recursion.  

It then easily follows that given the set partition $B_1, \ldots, B_k$
the sequences $(r_{i,1}, \ldots, r_{i,t_i})$, $i=1, \dots, k$ can be constructed in worst case time $O(n^2log(n))$ and that given the  
$(r_{i,1}, \ldots, r_{i,t_i})$, $i=1, \dots, k$, we can construct the set 
partition $B_1, \ldots, B_k$ in worst case time $O(n^2log(n))$. \hfil$\Box$

It then follows that our ranking procedure for 
$\vec{C}_{n,\vec{s}}$ requires $O(n^2log(n)$ bit operations.

The unranking procedure for $\vec{C}_{n\vec{s}}$ comes from simply reversing the ranking procedure using Theorem \ref{unrankdecfn} and Lemma \ref{produnrank}.
Again we will exhibit the procedure by finding the tree $T_1$ whose 
rank is 50,005 in $\vec{C}_{12,\langle 0,0,2,0,3,0,2,0,1,0,2,0
\rangle}$.   The first step is to carry out the series of quotients and 
remainders according to Lemma \ref{produnrank}. In our case, this leads to the following calculations. 
\begin{eqnarray*}
50,005 &=& (6 \times 7560) + 4645 \\
4645 &=& (22 \times 210) + 25\\
25 &=& (2 \times 10) + 5 \\
5 &=& (5 \times 1) + 0.
\end{eqnarray*}
It then follows that we can construct the sequence corresponding to 
$\overline{\pi}_{T_1}$ by concatonating the sequences 
$\vec{u}_1, \ldots, \vec{u}_5$ where 
\begin{enumerate}
\item $\vec{u}_1$ is the decreasing function of rank 6 in ${\cal DF}_{10,1}$, 
\item $\vec{u}_2$ is the decreasing function of rank 22 in ${\cal DF}_{9,2}$, 
\item $\vec{u}_3$ is the decreasing function of rank 2 in ${\cal DF}_{7,2}$, 
\item $\vec{u}_4$ is the decreasing function of rank 5 in ${\cal DF}_{5,2}$ and 
\item $\vec{u}_5$ is the decreasing function of rank 0 in ${\cal DF}_{3,3}$. 
\end{enumerate}

It is clear that the sequence of rank 6 in ${\cal DF}_{10,1}$ is $
\langle 7 \rangle$.

To find the element $\langle f(1),f(2) \rangle$ of rank 22 in ${\cal DF}_{9,2}$, we use the procedure 
in Theorem \ref{unrankdecfn}. We start by setting $m' := 22$ and $s = 2$. 
Since $\binom{7}{2} = 21 < 22 < \binom{8}{2} = 28$, then 
$f(1) -1 = 7$ and hence $f(1) =8$.  Then we set $m' := 22 -21 = 1$ and 
$s =1$.  Since $1 -\binom{1}{1} =0$, we get that $f(2) -1 =1$ or $f(2) =2$. 
Thus $\langle 8,2 \rangle$ has rank 22 in ${\cal DF}_{9,2}$.

To find the element $\langle f(1),f(2) \rangle$ of rank 2 in ${\cal DF}_{7,2}$, we again use the procedure 
in Theorem \ref{unrankdecfn}. We start by setting $m' := 2$ and $s = 2$. 
Since $\binom{2}{2} = 1 < 2 < \binom{3}{2} = 3$, then 
$f(1) -1 = 2$ and hence $f(1) =3$.  Then we set $m' := 2 -1 =1$ and 
$s =1$.  Since $1 -\binom{1}{1} =0$, we get that $f(2) -1 =1$ or $f(2) =2$. 
Thus $\langle 3,2 \rangle$ has rank 2 in ${\cal DF}_{7,2}$.

To find the element $\langle f(1),f(2) \rangle$ of rank 5 in ${\cal DF}_{5,2}$, we again use the procedure 
in Theorem \ref{unrankdecfn}. We start by setting $m' := 2$ and $s = 2$. 
Since $\binom{3}{2} = 3 < 5 < \binom{4}{2} = 6$,  
$f(1) -1 = 3$ and hence $f(1) =4$.  Then we set $m' := 5 -3 =2$ and 
$s =1$.  Since $2 -\binom{2}{1} =0$, we get that $f(2) -1 =2$ or $f(2) =3$. 
Thus $\langle 4,3 \rangle$ has rank 2 in ${\cal DF}_{5,2}$.

Finally there is only one element in ${\cal DF}_{3,3}$ which is 
$\langle 3,2,1\rangle$. Since the last step is alway trivial, it is most efficient to have the last sequence be as long as possible.  This is why we 
order the sizes of the set partition by increasing order.

Thus the sequence corresponding to the tree $\overline{\pi}_{T_1}$ is 
$$\langle 7, 8, 2, 3, 2, 4, 3, 3, 2, 1 \rangle.$$
It is easy to reconstruct $\overline{\pi}_{T_1}$ from this sequence and hence 
$$\overline{\pi}_{T_1} = \langle \{8\}, \{10,3\}, \{5,4\},\{9,7\},\{11,6,2\}\rangle.$$
It then follows that 
$$\pi_{T_1} = 
\langle \emptyset,\emptyset, \{10,3\}, \emptyset, \{11,6,2\}, \emptyset, 
\{5,4\}, \emptyset, \{8\}, \emptyset, \{9,7\}, \emptyset \rangle.$$
The function $f_1$ corresponding to $\pi_{T_1}$ and its image under 
$\Theta$ are pictured in Figure 6.

\fig{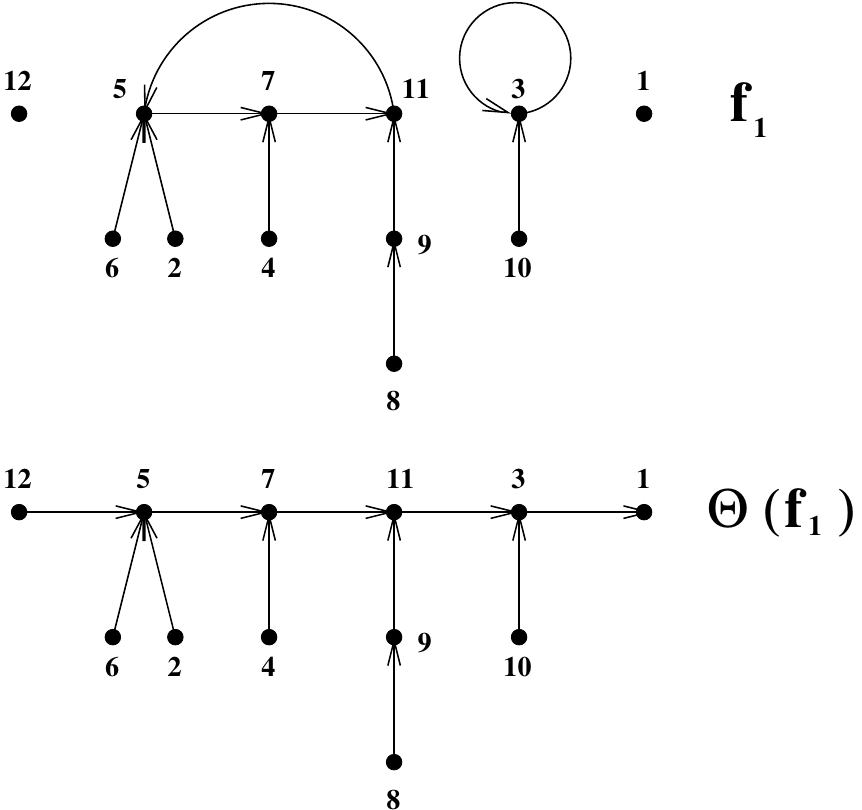}{The tree of rank 50,005 in $\vec{C}_{12,
\langle 0,0,2,0,3,0,2,0,1,0,2,0\rangle}$}

%\begin{figure}[htbp]
%\vspace{1.5in}
%\fcaption{Labeled tree {\it T.}}
%\end{figure}
%\begin{center}
%\
%\psfig{figure=T1.eps,width=2in,height=1.5in}
%\fcaption{The tree of rank 50,005 in $\vec{C}_{12,
%\langle 0,0,2,0,3,0,2,0,1,0,2,0\rangle}$}
%\end{center}
%\end{figure}

The problem of ranking and unranking trees with a given multiset of degrees is just an extension of the the problem ranking and unranking trees with a given sequence of degrees.  That is, the distribution of degrees is just another set partition.  For example, consider the sequence of degrees for the tree $T_0$ pictured in Figure 5, $\vec{s} = \langle 0, 0, 2, 0, 3, 0, 2, 0, 1, 0, 2, 0, 0 \rangle$. We can view that sequence as a set partition, $\Delta(\vec{s}) = \langle \Delta_0, \ldots, \Delta_3\rangle$ where $\Delta_i$ 
is the set of places where 
$i$ appears in the sequence $\vec{s}$. In our example, we would 
identify $\vec{s}$ with the set partition 
$$ \Delta(\vec{s}) =\langle \{1,2,4,6,8,10,12\}, \{9\}, \{3,7,11\}, \{5\}\rangle.$$

Just as in the case where we ranked and unranked set partitions associated 
with trees in $\vec{C}_{n,\vec{s}}$, it is more efficient if we rearrange the set partition by increasing size.  This means that we must have a data structure to record the degrees associated with the set partition which in this case is just the triples $(\Delta_i, |\Delta_i|, i)$. It is easy to see that we can produce such a list in linear time from the tree.  In our example, 
we would produce the list
$$ \Delta(\vec{s}) =\langle (\{1,2,4,6,8,10,12\},7,0), (\{9\},1,1), 
(\{3,7,11\},3,2), (\{5\},1,3)\rangle.$$
Using a lexicographic bucket sort algorithm \cite{W}, we can sort this list according to the lexicographic order on the last two entries of the triples to produce the list 
$$\overline{\Delta(\vec{s})} = \langle (\{9\},1,1), (\{5\},1,3),  
(\{3,7,11\},3,2),(\{1,2,4,6,8,10,12\},7,0) \rangle.$$
Then we use this ordering to produce an ordered set partition 
$\pi_{\vec{s}}$ where we ignore any empty partitions. In our example, 
we would produce 
$$\pi_{\vec{s}} = 
\{9\}, \{5\}, \{3,7,11\},\{1,2,4,6,8,10,12\}\rangle.$$
Finally, we use this set partition to produce a sequence of 
decreasing functions. That is, we let $E=\{1, \ldots , n\}$. 
Let $\overline{\pi} = (A_1, \ldots, A_k)$ be an ordered
set partition of $E$ where each block $E_i$ is nonempty and ordered in
decreasing order, $A_i=\{a_{i,1}> \ldots > a_{i,t_i}\}$, that comes from some 
element $\pi$ in $\pi_{\vec{s}}$ as described above.
Let $(r_{1,1}, \ldots, r_{1,t_1})$ be the ordered sequence of ranks of
the respective $(a_{1,1}, \ldots, a_{1,t_1})$ in $E$.
In general, let $(r_{i,1}, \ldots, r_{i,t_i})$ be the ranks of the
respective $(a_{i,1}, \ldots, a_{i,t_i})$ in
$E\setminus \cup_{j=1}^{i-1}A_j$.
For our example, \\
$(r_{1,1}) =  (9)$ can be considered an element of 
${\cal DF}_{12,1}$\\
$(r_{2,1}) = (5)$ can be considered an element of 
${\cal DF}_{11,1}$\\
$(r_{3,1},r_{3,2},r_{3,3}) = (9,6,3)$ can be considered an element of 
${\cal DF}_{10,3}$\\
$(r_{4,1},r_{4,2}, r_{4,3},r_{4,4},r_{4,5},r_{4,6},r_{4,7}) = 
(7,6,5,4,3,2,1)$ can be considered an element of 
${\cal DF}_{7,7,2}$.\\ 
Thus we produce the sequence 
$\overline{\pi}_{\vec{s}} = \langle 9,5,9,6,3,,7,6,5,4,3,2,1\rangle
$ to code the sequence $\vec{s}$ which can be considered an element 
of ${\cal DF}_{12,1}\times {\cal DF}_{11,1}\times {\cal DF}_{10,3} \times {\cal DF}_{7,7}$.  We can then concatenate this sequence $\overline{\pi}_{\vec{s}}$ with the sequence $\overline{\pi}_{T_0}$ to produce a sequence 
$\Pi_{\vec{s},T_0}$.  In our example,  
$$\Pi_{\vec{s},T_0}= \langle 9,5,9,6,3,,7,6,5,4,3,2,1 6,8,5,6,2,3,1,3,2,1 \rangle$$ coming from an element 
of 
$${\cal DF}_{12,1}\times {\cal DF}_{11,1}\times {\cal DF}_{10,3} \times {\cal DF}_{7,7} \times {\cal DF}_{10,1}\times {\cal DF}_{9,2} \times {\cal DF}_{7,2} \times 
{\cal DF}_{5,2}\times {\cal DF}_{3,3}$$ or as a leaf in the tree  
$$T_{{\cal DF}_{12,1}}\otimes T_{{\cal DF}_{11,1}}\otimes T_{{\cal DF}_{10,3}} 
\otimes T_{{\cal DF}_{7,7}} \otimes 
T_{{\cal DF}_{10,1}}\otimes T_{{\cal DF}_{9,2}} \otimes T_{{\cal DF}_{7,2}} 
\otimes T_{{\cal DF}_{5,2}} \otimes T_{{\cal DF}_{3,3}}.$$
Note that the size of the trees needed for the product lemma, Lemma 
\ref{prodrank}, are \\
$|T_{{\cal DF}_{3,3}}| = \binom{3}{3} = 1$,\\
$|T_{{\cal DF}_{5,2}} \otimes T_{{\cal DF}_{3,3}}| = \binom{5}{2} \cdot \binom{3}{3} = 10$,\\
$|T_{{\cal DF}_{7,2}} 
\otimes T_{{\cal DF}_{5,2}} \otimes T_{{\cal DF}_{3,3}}| = 
\binom{7}{2} \cdot \binom{5}{2} \cdot \binom{3}{3} = 210$,\\
$|T_{{\cal DF}_{9,2}} \otimes T_{{\cal DF}_{7,2}} 
\otimes T_{{\cal DF}_{5,2}} \otimes T_{{\cal DF}_{3,3}}
T_{{\cal DF}_{7,2}} 
\otimes T_{{\cal DF}_{5,2}} \otimes T_{{\cal DF}_{3,3}}| =$\\ 
$\binom{9}{2} \cdot \binom{7}{2} \cdot \binom{5}{2} \cdot \binom{3}{3} = 
7560$,\\
$|T_{{\cal DF}_{10,1}}\otimes T_{{\cal DF}_{9,2}} \otimes T_{{\cal DF}_{7,2}} 
\otimes T_{{\cal DF}_{5,2}} \otimes T_{{\cal DF}_{3,3}}| =$\\
$\binom{10}{1} \cdot \binom{9}{2} \cdot \binom{7}{2} \cdot \binom{5}{2} \cdot \binom{3}{3} = 
75600$,\\
$|T_{{\cal DF}_{7,7}} \otimes T_{{\cal DF}_{10,1}}\otimes T_{{\cal DF}_{9,2}} \otimes T_{{\cal DF}_{7,2}} 
\otimes T_{{\cal DF}_{5,2}} \otimes T_{{\cal DF}_{3,3}}| =$\\
$\binom{7}{7} \cdot \binom{10}{1} \binom{9}{2} \cdot \binom{7}{2} \cdot \binom{5}{2} \cdot \binom{3}{3} = 
75600$,\\
$|T_{{\cal DF}_{10,3}} \otimes T_{{\cal DF}_{7,7}} \otimes T_{{\cal DF}_{10,1}}\otimes T_{{\cal DF}_{9,2}} \otimes T_{{\cal DF}_{7,2}} 
\otimes T_{{\cal DF}_{5,2}} \otimes T_{{\cal DF}_{3,3}}| =$ \\
$\binom{10}{3} \cdot 
\binom{7}{7} \cdot \binom{10}{1} \binom{9}{2} \cdot 
\binom{7}{2} \cdot \binom{5}{2} \cdot \binom{3}{3} = 
9,072,000$,\\
$|T_{{\cal DF}_{11,1}} 
\otimes T_{{\cal DF}_{10,3}} 
\otimes T_{{\cal DF}_{7,7}} \otimes T_{{\cal DF}_{10,1}}\otimes T_{{\cal DF}_{9,2}} \otimes T_{{\cal DF}_{7,2}} 
\otimes T_{{\cal DF}_{5,2}} \otimes T_{{\cal DF}_{3,3}}| =$ \\
$\binom{11}{1} \cdot \binom{10}{3} \cdot 
\binom{7}{7} \cdot \binom{10}{1} \binom{9}{2} \cdot \binom{7}{2} \cdot \binom{5}{2} \cdot \binom{3}{3} = 
99,792,000$,\\
$|T_{{\cal DF}_{12,1}} \otimes T_{{\cal DF}_{11,1}} 
\otimes T_{{\cal DF}_{10,3}} 
\otimes T_{{\cal DF}_{7,7}} \otimes T_{{\cal DF}_{10,1}}\otimes T_{{\cal DF}_{9,2}} \otimes T_{{\cal DF}_{7,2}} 
\otimes T_{{\cal DF}_{5,2}} \otimes T_{{\cal DF}_{3,3}}| =$ \\
$\binom{12}{1} \cdot \binom{11}{1} \cdot \binom{10}{3} \cdot 
\binom{7}{7} \cdot \binom{10}{1} \binom{9}{2} \cdot \binom{7}{2} \cdot \binom{5}{2} \cdot \binom{3}{3} = 
1,197,504,000$.\\
Thus there are a total of 1,197,504,000 trees in $\vec{C}_{12,1}$ whose 
degree sequence yields the multiset $S= (0^7,1^1,2^3,3^1)$. We can then use 
the product lemma, Lemma \ref{prodrank}, to compute the rank of 
$T_0$ in $\vec{C}_{12,S}$ as follows. 
\begin{eqnarray*}
rank_{\vec{C}_{12,S}}(T_0) &=& rank_{{\cal DF}_{12,1}}(\langle 9 \rangle) \times 99,720,000\\
&&+ rank_{{\cal DF}_{11,1}}(\langle 5 \rangle) \times 9,072,000\\
&&+ rank_{{\cal DF}_{10,3}}(\langle 9,6,3 \rangle) \times 75,600\\
&&+ rank_{{\cal DF}_{7,7}}(\langle 7,6,5,4,3,2,1 \rangle) \times 75,600\\
&&+ rank_{{\cal DF}_{10,1}}(\langle 5 \rangle) \times 7,560\\
&&+ rank_{{\cal DF}_{9,2}}(\langle 8,5 \rangle) \times 210\\
&&+ rank_{{\cal DF}_{7,2}}(\langle 6,2 \rangle) \times 10\\
&&+ rank_{{\cal DF}_{5,2}}(\langle 3,1 \rangle) \times 1.
\end{eqnarray*}
By Theorem \ref{rankdecfn}, we have that \\
$rank_{{\cal DF}_{12,1}}(\langle 9 \rangle) = \binom{9-1}{1} = 8$,\\
$rank_{{\cal DF}_{11,1}}(\langle 5 \rangle) = \binom{5-1}{1} = 4$,\\
$rank_{{\cal DF}_{10,3}}(\langle 9,6,3 \rangle) = 
\binom{9-1}{3} + \binom{6-1}{2} + \binom{3-1}{1} = 56+10+2=68$,\\
$rank_{{\cal DF}_{7,7}}(\langle7,6,5,4,3,2,1 \rangle) = 0$,\\
$rank_{{\cal DF}_{10,1}}(\langle 5 \rangle) = \binom{5-1}{1} = 4$,\\
$rank_{{\cal DF}_{9,2}}(\langle 8,5 \rangle) = \binom{8-1}{2} + 
\binom{5-1}{1} = = 21 +4 = 25$,\\
$rank_{{\cal DF}_{7,2}}(\langle 6,2 \rangle) =  \binom{6-1}{2} + \binom{2-1}{1} = 10 +1 = 11$,\\
$rank_{{\cal DF}_{5,2}}(\langle 3,1 \rangle)= \binom{3-1}{2} + \binom{1-1}{1} = 1 +0 = 1$.\\
Thus 
\begin{eqnarray*}
rank(\overline{\pi}_0) &=& (9 \times 99,720,000) + 
(5 \times 9,072,00) + (68 \times 75,600) + 
(0 \times 75,600) \\
&&+(4\times 7560) +(25 \times 210)+ (11\times 10) + 
(1 \times 1) = 843,342,641.
\end{eqnarray*}
Thus the tree $T_0$ pictured in Figure 5 has rank 843,342,641 among all 
the trees $T \in \vec{C}_{12,(0^7,1^1,2^3,3^1)}$.

The unranking procedure for $\vec{C}_{n,S}$ comes from simply reversing the ranking procedure using Theorem \ref{unrankdecfn} and Lemma \ref{produnrank}.
Again we will exhibit the procedure by finding the tree $T_2$ whose 
rank is 60,000,00  in $\vec{C}_{12,(0^7,1^1,2^3,3^1)}$.   The first step is to carry out the series of quotients and 
remainder according to Lemma \ref{produnrank}. In our case, this leads to the following calculations. 
\begin{eqnarray*}
60,000,00 &=& (6 \times 99,792,000) + 1,248,000\\
1,248,000 &=& (0 \times 9,072,000) + 1,248,000\\
1,248,000 &=& (16 \times 75,600) + 38,400\\
38,400 &=& (0 \times 75,600) + 38,400\\
38,400 &=& (5 \times 7560) + 600 \\
600 &=& (2 \times 210) + 180\\
180 &=& (18 \times 10) + 0 \\
0 &=& (0 \times 1) + 0.
\end{eqnarray*}
It then follows that we can construct the sequence corresponding to 
$\overline{\pi}_{,\vec{s},T_2}$ by concatonating the sequences 
$\vec{v}_1, \ldots, \vec{v}_9$ where 
\begin{enumerate}
\item $\vec{v}_1$ is the decreasing function of rank 6 in ${\cal DF}_{12,1}$, 
\item $\vec{v}_2$ is the decreasing function of rank 0 in ${\cal DF}_{11,1}$, 
\item $\vec{v}_3$ is the decreasing function of rank 16 in ${\cal DF}_{10,3}$, 
\item $\vec{v}_4$ is the decreasing function of rank 0 in ${\cal DF}_{7,7}$ and \item $\vec{v}_5$ is the decreasing function of rank 5 in ${\cal DF}_{10,1}$.
\item $\vec{v}_6$ is the decreasing function of rank 2 in ${\cal DF}_{9,2}$.
\item $\vec{v}_7$ is the decreasing function of rank 18 in ${\cal DF}_{7,2}$.
\item $\vec{v}_8$ is the decreasing function of rank 0 in ${\cal DF}_{5,2}$.
\item $\vec{v}_9$ is the decreasing function of rank 0 in ${\cal DF}_{3,3}$.

\end{enumerate}
\ \\
It is clear that the sequence of rank 6 in ${\cal DF}_{12,1}$ is $
\langle 7 \rangle$ and the sequence of rank 0 in ${\cal DF}_{11,1}$ is $
\langle 1 \rangle$.\\
\ \\
To find the element $\langle f(1),f(2),f(3) \rangle$ of rank 16 in 
${\cal DF}_{10,3}$, we use the procedure 
in Theorem \ref{unrankdecfn}. We start by setting $m' := 16$ and $s = 3$. 
Since $\binom{5}{3} = 10 < 16 < \binom{6}{3} = 20$, then 
$f(1) -1 = 5$ and hence $f(1) =6$.  
Then we set $m' := 16 -10 = 1$ and 
$s =2$.  Since $\binom{4}{2} = 6 \leq 6 < \binom{5}{2} = 10$, then 
$f(2) -1 = 4$ and hence $f(2) =5$. Finally we set $m' := 6 -6 = 0$ and 
$s =1$. Since $0 -\binom{0}{1} =0$, we get that $f(3) -1 =0$ or $f(3) =1$. 
Thus $\langle 6,5,1 \rangle$ has rank 16 in ${\cal DF}_{10,3}$.\\
\ \\
It is clear that the sequence of rank 0 in ${\cal DF}_{7,7}$ is $
\langle 7,6,5,4,3,2,1 \rangle$. \\
\ \\
Next it is clear that element of rank 5 in ${\cal DF}_{10,1}$ is 
$\langle 6 \rangle$.\\
\ \\
To find the element $\langle f(1),f(2) \rangle$ of rank 2 in ${\cal DF}_{9,2}$, we use the procedure 
in Theorem \ref{unrankdecfn}. We start by setting $m' := 2$ and $s = 2$. 
Since $\binom{2}{2} = 1 < 2 < \binom{3}{2} = 3$, then 
$f(1) -1 = 2$ and hence $f(1) =3$.  Then we set $m' := 2 -1 = 1$ and 
$s =1$.  Since $1 -\binom{1}{1} =0$, we get that $f(2) -1 =1$ or $f(2) =2$. 
Thus $\langle 3,2 \rangle$ has rank 2 in ${\cal DF}_{9,2}$.\\
\ \\
To find the element $\langle f(1),f(2) \rangle$ of rank 18 in ${\cal DF}_{7,2}$, we again use the procedure 
in Theorem \ref{unrankdecfn}. We start by setting $m' := 18$ and $s = 2$. 
Since $\binom{6}{2} = 15 < 18 < \binom{7}{2} = 21$, then 
$f(1) -1 = 6$ and hence $f(1) =7$.  Then we set $m' := 18-15 =3$ and 
$s =1$.  Since $3 -\binom{3}{1} =0$, we get that $f(2) -1 =3$ or $f(2) =4$. 
Thus $\langle 7,4 \rangle$ has rank 18 in ${\cal DF}_{7,2}$.\\
\ \\
Finally the sequence of rank 0 in ${\cal DF}_{7,2}$ is clearly $\langle 2,1 \rangle$ and the element of rank 0 in ${\cal DF}_{3,3}$ is $\langle 3,2,1 \rangle$.\\
\ \\
Thus the sequences corresponding to the set partitions 
$\overline{\pi}_{\vec{s}}$ 
and $\overline{\pi}_{T_2}$ are 
\begin{eqnarray*}
\overline{\pi}_{\vec{s}}&:&\langle 7,1,6,5,1,7,6,5,4,3,2,1\rangle \ 
\mbox{and}\\
\overline{\pi}_{T_2}&:&\langle 6,3,2,7,4,2,1,3,2,1\rangle.
\end{eqnarray*}

It is easy to reconstruct $\overline{\pi}_{S}$ 
and $\overline{\pi}_{T_2}$ to get that 
\begin{eqnarray*}
\overline{\pi}_{\vec{s}} &=& \langle \{7\}, \{1\}, \{8,6,2\},\{12,11,10,9,5,4,3\}\rangle \ \mbox{and}\\
\overline{\pi}_{T_1} &=& \langle \{7\}, \{4,3\}, \{11,8\},\{5,2\},\{10,9,6\}\rangle
\end{eqnarray*}

It then follows that
$$\vec{s} = \langle 3,2,0,0,0,2,1,2,0,0,0,0\rangle$$
and  
$$\pi_{T_1} = 
\langle \{10,9,6\},\{4,3\}, \emptyset, \emptyset, \emptyset, 
\{11,8\}, \emptyset, \{7\}, \emptyset, \{5,2\}, \emptyset,\emptyset,\emptyset  \rangle.$$
Thus the function $f_2$ corresponding to $\pi_{T_3}$ and its image under 
$\Theta$ are pictured in Figure~7.

\fig{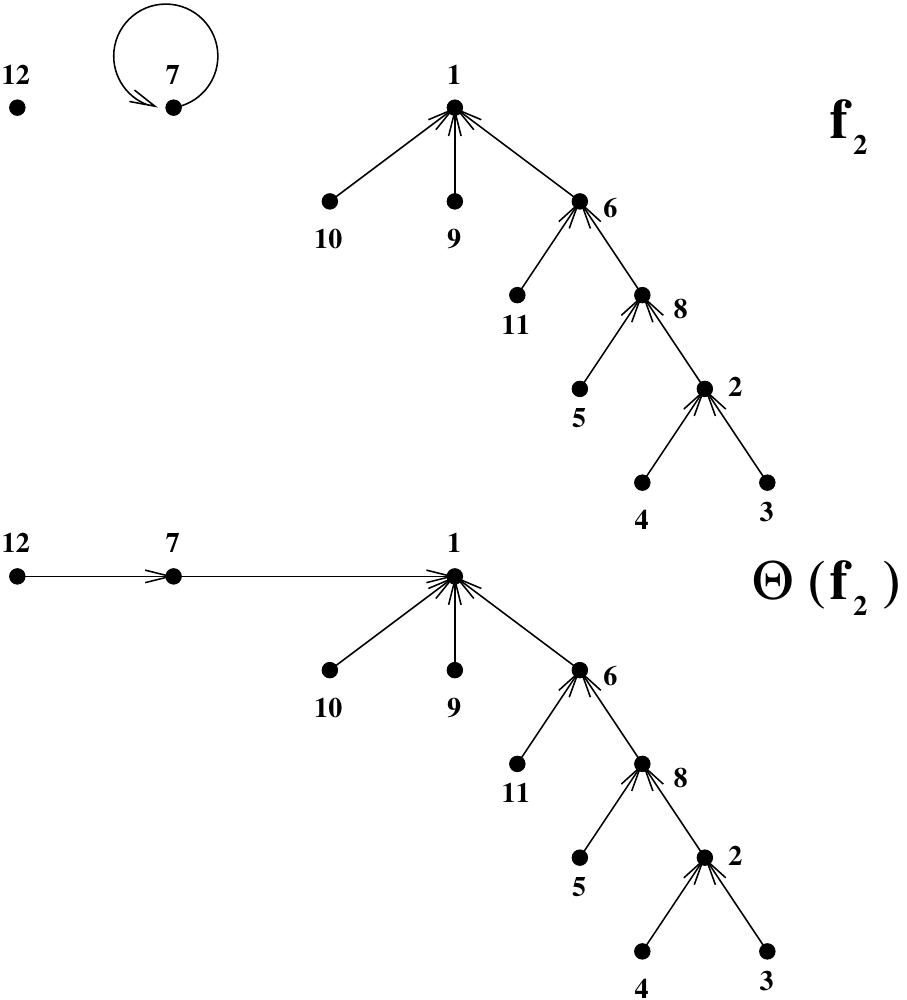}{The tree of rank 60,000,000 in $\vec{C}_{12,
(0^7,1^1,2^3,3^1)}$}

%\begin{figure}[htbp]
%\vspace{1.5in}
%\fcaption{Labeled tree {\it T.}}
%\end{figure}
%\begin{center}
%\
%\psfig{figure=T2.eps,width=2in,height=2in}
%\fcaption{The tree of rank 60,000,000 in $\vec{C}_{12,
%(0^7,1^1,2^3,3^1)}$}
%\end{center}
%\end{figure}

We note that essentially the same analysis of the complexity of 
ranking and unranking relative $\vec{C}_{n,\vec{s}}$ applies to 
the complexity of ranking and unranking relative to $\vec{C}_{n,S}$ for 
any multiset $S$ so that it requires $O(n^2log(n))$ bit operations to 
rank and unrank relative $\vec{C}_{n,S}$.

%\nonumsection{References}
%References are to be listed in the order cited in the text. Use the style
%shown in the following examples. For journal names, use the standard
%abbreviations. Typeset references in 9 pt Times Roman.

\end{document}